\def\marker{\>\hbox{${\vcenter{\vbox{
    \hrule height 0.4pt\hbox{\vrule width 0.4pt height 6pt
    \kern6pt\vrule width 0.4pt}\hrule height 0.4pt}}}$}\>}

\documentclass[12pt]{article}

\usepackage[left=2cm,right=2cm,top=3cm,bottom=3cm]{geometry}

\usepackage{amsmath, amssymb, latexsym, amsthm}

\usepackage{graphicx}
\usepackage{fullpage}

\usepackage{tikz}
\usetikzlibrary{shapes}

\newtheorem{theorem}{Theorem} 
\newtheorem{theorem*}{Theorem} 

\newtheorem{corollary}[theorem]{Corollary}

\theoremstyle{definition}

\theoremstyle{remark}

{\end{enumerate}}

\usepackage[hang,flushmargin]{footmisc}

\title{ 3-colorability of graphs with minimum degree at least 6}
\author{Nicholas Crawford\footnotemark[1] and Sogol Jahanbekam\footnotemark[2]}

\date{}

\begin{document}

\maketitle

\begin{abstract}
Let $G$ be an $n$-vertex graph and  let $L:V(G)\rightarrow P(\{1,2,3\})$ be a list assignment over the vertices of $G$, where each vertex with list of size 3 and of degree at most 5 has at least three neighbors with lists of size 2. We can determine $L$-choosability of $G$ in $O(1.3196^{n_3+.5n_2})$ time, where $n_i$ is the number of vertices in $G$ with list of size $i$ for $i\in \{2,3\}$. As a corollary, we conclude that the 3-colorability of any graph $G$ with minimum degree at least 6 can be determined in $O(1.3196^{n-.5\Delta(G)})$ time.

\end{abstract}

\renewcommand{\thefootnote}{\fnsymbol{footnote}}
\footnotetext[2]{
Department of Mathematics and Statistics, San Jose State University, San Jose, CA; {\tt sogol.jahanbekam@sjsu.edu.}\\{Research was supported in part by by NSF grant CMMI-1727743.}
}

\footnotetext[2]{
 Department of Mathematics and Statistics, San Jose State University, San Jose, CA;
{\tt Nicholas.crawford@sjsu.edu }. }
\renewcommand{\thefootnote}{\arabic{footnote}}

\baselineskip18pt

\section{Introduction}
Let  $G$ be a graph. We denote the vertex set of $G$ by $V(G)$, the edge set of $G$ by $E(G)$, and  the minimum degree of $G$ as $\delta(G)$.
For a vertex $v$ in $G$, the \textit{open neighborhood} of $v$ in $G$, denoted by $N(v)$, is the set of neighbors of $v$ in $G$. 

 We say that a coloring of vertices of $G$ is \emph{proper} if all adjacent vertices receive different colors. A graph $G$ is \emph{k-colorable} if it has a proper coloring using at most \emph{k} colors. The \emph{chromatic number} of a graph G, denoted $\chi(G)$, is the smallest integer \emph{k} such that G is \emph{k}-colorable.

Suppose for each vertex \emph{v} in $V(G)$, there exist a list of colors denoted $L(v)$. A \emph{proper list coloring} of G is a choice function that maps every vertex \emph{v} to a color in $L(v)$ in such a way that the coloring is proper. A graph is $k$-choosable if it has a proper list coloring whenever each vertex has a list of size $k$.  For additional definitions we refer the reader to \cite{W}.


By the definition, a graph $G$ is 1-colorable if and only if it has no edge. We can also apply a simple Breadth First Search algorithm to determine 2-coloability of a graph in polynomial time \cite{W}. The $k$-colorability problem for $k\geq 3$ is an NP-complete problem \cite{MTB}, and as a result determining the chromatic number of a graph is NP-complete.

In 1971 Christofides discovered the first non-trivial algorithm for computing the chormatic number of an $n$-vertex graph with complexity $n!n^{O(1)}$ \cite{C}. The best known result for determining the chromatic number of a graph was done by Bj\"orklund, Husfeldt, and Koivisto in 2009 with complexity $2^{n}O(1)$ \cite{BHK}. 

The $k$-coloring problem  is a highly studied problem in graph theory. Improving the complexity for $k$-coloring problem, even for small values of $k$ like 3 and 4,  could lead to improved complexity of the general chromatic number problem.  In 2013,  Eppstein and Beigal \cite{BE} created an algorithm that determined whether a graph was 3-colorable in $O(1.3289^{n})$ time.  This algorithm is the lowest complexity algorithm to date for 3-colorability of general graphs. Crawford et. al. proved in \cite{CJP}  that you can determine the 3 colorability of a graph with minimum degree at least 7 and 8 in $O(1.32^{n-.73\Delta(G)})$ and $O(1.3158^{n-.7\Delta(G)})$, respectively.


 In this paper, we prove the following results.

 \begin{theorem}\label{main}
Let $G$ be an $n$-vertex graph with minimum degree at least 6.  Let $L:V(G)\rightarrow P{\{1,2,3\}}$ be a list assignment over the vertices of $G$. We can determine $L$-choosability of $G$ in $O(1.3196^{n_3+.5n_2})$ time, where $n_i$ is the number of vertices in $G$ with list of size $i$ for $i\in \{2,3\}$.
\end{theorem}

Note that if a vertex in a proper coloring of a graph $G$ using three colors $1,2,3$ gets color 1, then each of its neighbors must get a color in $\{2,3\}$. This observation and  Theorem \ref{main} imply the following corollary.

\begin{corollary}
Let $G$ be an $n$-vertex graph with  minimum degree at least $6$. We  determine in $O(1.3196^{n-.5\Delta(G)})$ time  if $G$ is 3-colorable. 
\end{corollary}

To prove Theorem \ref{main} we prove the following stronger theorem.

 \begin{theorem}\label{stronger}
Let $G$ be an $n$-vertex graph and  let $L:V(G)\rightarrow P(\{1,2,3\})$ be a list assignment over the vertices of $G$. If each vertex of degree at most 5 in $G$, with a list of size 3, has at least three neighbors with lists of size 2, then we can determine $L$-choosability of $G$ in $O(1.3196^{n_3+.5n_2})$ time, where $n_i$ is the number of vertices with list of size $i$ in $G$.
\end{theorem}

\section{Proof of Theorem \ref{stronger}}

We apply induction on $k$, the number of vertices of the graph, to prove that if $L:V(G)\rightarrow P(\{1,2,3\})$ is a list assignment over the vertices of a $k$-vertex graph $G$  with the condition that  each vertex of degree at most 5 in $G$, with a list of size 3, has at least three neighbors with lists of size 2, then we can determine $L$-choosability of $G$ in  $O(1.3196^{k_3+.5k_2})$ time, where $k_i$ is the number of vertices in $G$ with lists of size $i$.

If $k=1$, then $G=K_1$, which is $L$-choosable. Therefore in this case we need no operation to determine $L$-choosability of $G$, where $L$ gives at least 1 and  at most 3 colors from $\{1,2,3\}$ to each vertex.

Now suppose the assertion holds for any graph satisfying the condition of the problem and having smaller than $n$ vertices. Let $G$ be an $n$-vertex graph. Suppose \\ $L:V(G)\rightarrow P(\{1,2,3\})$ is a list assignment over the vertices of $G$ in such a way that each vertex of degree at most 5, and list of size 3, has at least 3 neighbors with lists of size 2. 

 If $G$ has a vertex $u$ with a list of size 1, then we can remove the color in $L(u)$ from the list of colors in the neighborhood of $u$ and then study the smaller graph $G-u$. Hence we may suppose $G$ has no vertex with a list of size 1. Suppose $G$ has $n_3$ vertices with lists of size 3 and  $n_2$ vertices with lists of size 2.

We consider four cases. Note that each of the cases below is considering a specific structure in graph $G$, which can be found via a polynomial-time algorithm.

\noindent 
\textbf{Case 1}: A vertex $v$ with list size $3$ in $G$ is adjacent to at least 6 vertices with list size 3 (Figure 1).

\begin{figure}[htp]
    \centering
    \includegraphics[width=7cm]{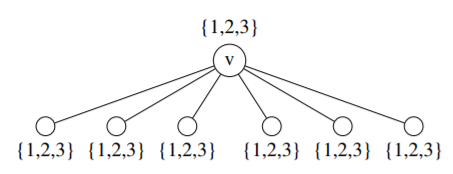}
    \caption{Graph $G$ with List $L$.}
    \label{TL1}
\end{figure}

There are three possibilities for the color of $v$: 1, 2, or 3. First suppose $v$ gets color 1 in $G$. In this case remove color 1 from the lists on the neighbors of $v$ and let $L'$ be the resulting list assignment (Figure 2). In this case $G$ is $L$-choosable with vertex $v$ getting color 1 if and only if $G-v$ is $L'$-choosable.

\begin{figure}[htp]
    \centering
    \includegraphics[width=7cm]{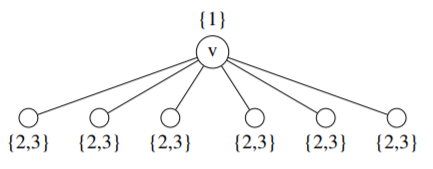}
    \caption{When vertex $v$ gets color 1.}
    \label{TL2}
\end{figure}

Note that at least 6 vertices in the neighborhood of $v$ have lists of size 3 in $L$ but lists of size 2 in $L'$, because we have removed color 1 from their lists. Moreover, $v$ has a list of size 3 in $G$, while it is not a vertex in $G-v$. Therefore $G-v$ contains at most $n_3-7$ vertices with lists of size 3. Since $G-v$ has $n-1$ vertices, by the induction hypothesis we can determine in  $O((1.3196)^{n_3-7+.5(n_2+6)})$ time if it is $L'$-choosable.

A similar argument as above applies when $v$ gets color 2 or color 3. Hence we can determine in  $O(3(1.3196)^{n_3-7+.5(n_2+6)})$ time if $G$ is $L$-choosable. Since $3(1.3196)^{n_3-7+.5(n_2+6)}< (1.3196)^{n_3+.5n_2}$, the assertion holds in this case.

\noindent
\textbf{Case 2}: A vertex $v$ with list size 3 is adjacent to $5$ vertices with lists of size 3.

By Case 1 and by the structure of $G$ we may suppose that $v$ has at least one  neighbor with list of size 2. Let $u_1,\ldots,u_5$ be neighbors of $v$ with lists of size 3 and let $u_6$ be a neighbor of $v$ with list of size 2. By symmetry we may suppose $L(u_6)=\{1,2\}$. We consider three subcases.

First suppose $u_6$ has another neighbor $w$ with a list containing 2 (i.e. $L(w)=\{1,2\}$, or $L(w)=\{2,3\}$, or $L(w)=\{1,2,3\}$). 
Note that the vertex $v$ can get color 1, 2, or 3.  
If $v$ gets color 1 in $G$, then the neighbors of $v$ cannot be color 1.
Define $L'(u)=L(u)-\{1\}$ for all neighbors $u$ of $v$. Hence $L(u_6)=\{2\}$. 
Therefore $u_6$ has to get color 2 in this case. As a result, $w$ cannot receive color $2$ anymore. Hence let $L'(w)=L(w)-\{2\}$. Thus $G$ is $L$-choosable with $v$ getting color $1$ if and only if $G-\{v,u_6\}$ is $L'$-choosable. Note that $G-\{v,u_6\}$ has $n-2$ vertices. Moreover, it contains either $n_3-7$ and $n_2+4$ vertices with lists of size 3 and 2, respectively, or  it contains  $n_3-6$ and $n_2+3$ vertices with lists of size 3 and 2, respectively (depending on the list size on vertex $w$).

 We can repeat this argument for the case $v$ gets color 2 and 3. If $v$ gets color $2$, then we can reduce the problem into $L'$-choosability of the graph $G-\{v,u_6\}$ with $n-2$ vertices, $n_3-6-i$ vertices with lists of size 3 and $n_2+5-j$ vertices with lists of size 2, where $\{i,j\}=\{0,1\}$. And when $v$ gets color $3$,  then we can reduce the problem into $L'$-choosability of the graph $G-\{v\}$ with $n-1$ vertices, $n_3-6$ vertices with lists of size 3 and $n_2+5$ vertices with lists of size 2. Hence we can determine in  $O((1.3196)^{n_3-6+.5(n_2+3)}+(1.3196)^{n_3-6+.5(n_2+4)}+(1.3196)^{n_3-6+.5(n_2+5)})$
time if $G$ is $L$-choosable. Since $(1.3196)^{-6+.5(3)}+(1.3196)^{-6+.5(4)}+(1.3196)^{-6+.5(5)}\leq 1$, the assertion holds in this case.

If $u_6$ has another neighbor $w$ with a list containing 1, we can repeat the above argument by switching the arguments for the case $v$ gets colors 1 and 2 to get the desired conclusion. Hence the final subcase is when $u_6$ has no other neighbor containing color $1$ or $2$. This implies that $u_6$ has degree 1 in $G$. As a result, $L$-choosability of $G$ can be reduced to list choosability of $G-\{v,u_6\}$, for three different list assignment corresponding to the cases that $v$ gets color 1, 2, or 3. At each of these cases, $G-\{v,u_6\}$ has $n_3-6$ vertices with lists of size 3 and $n_2+4$ vertices with lists of size 2. Therefore $L$-choosability of $G$ can be determined in time $O(3(1.3196)^{n_3-6+.5(n_2+4)})$ which is a subset of $O(1.3196^{n_3+.5n_2})$, as desired.

\noindent
\textbf{Case 3}: A vertex $v$ with list size 3 is adjacent to  4 vertices of lists of size 3. We need to consider more subcases here compared to Case 1 and Case 2. However, the idea of the proof is very similar to them. Therefore, to avoid redundancy we skip some details.

At each subcase we consider the three possibilities on the color of $v$. When $v$ gets color $i$ for $i\in \{1,2,3\}$, we reduce $L$-choosability of $G$ to an $L'$-choosability of a subgraph of $G$, which we call $G_i$. In the following subcases we explain what $G_i$ is in each case, but we skip the definition of $L'$, as $L'$ is naturally obtained from $L$ after removing the colors of the vertices with fixed colors from the list of  their neighbors.

By the structure of $G$ and since we may suppose Cases 1 and 2 do not happen,  the vertex $v$ has at least two neighbors with lists of size 2.  Let $u_1$ and $u_2$ be two such neighbors of $v$. By symmetry suppose the list of $u_1$ is $\{1,2\}$. Let the list on $u_2$ be $\{c_1,c_2\}$.

\noindent
\textbf{Subcase 1}: $u_1$ has another neighbor $u'_1$ with list of size 3 (Figure 3).

\begin{figure}[htp]
    \centering
    \includegraphics[width=7cm]{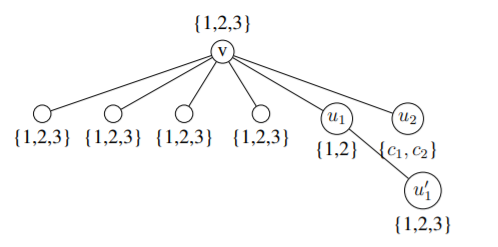}
    \caption{When $u_1$ has another neighbor with list $\{1,2,3\}$.}
    \label{TL3}
\end{figure}

Note that for the case that the color of $v$ is 1 or $2$, then we get only one possible color for $u_1$, and as a result the list of $u'_1$ decreases from size 3 to size 2 (Figure 4).

\begin{figure}[htp]
    \centering
    \includegraphics[width=7cm]{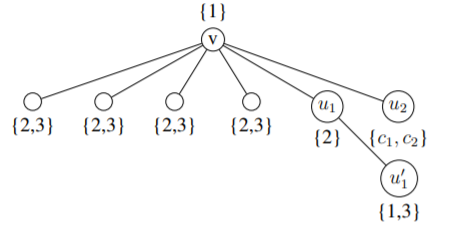}
    \caption{When $v$ gets color 1 in Subcase 1.}
    \label{TL4}
\end{figure}
 
  If $\{c_1,c_2\}=\{1,2\}$ and $u_2$ has degree 1 in $G$, then let $G_1=G_2=G-\{v,u_1,u_2\}$ and $G_3=G-\{v,u_2\}$.   Then  by the induction hypothesis we can determine $L$-choosability of $G$ in  
$O((1.3196)^{n_3-6+.5(n_2+3)}+(1.3196)^{n_3-6+.5(n_2+3)}+(1.3196)^{n_3-5+.5(n_2+3)})$ time, which are both subsets of $O(1.3196^{n_3+.5n_2})$.

   If $\{c_1,c_2\}=\{1,2\}$ and $u_2$ has degree at least 2 in $G$, then we may suppose $u_2$ has a neighbor with a list containing 2 (Otherwise it must contain a neighbor with a list containing 1, which implies a similar proof). In this case, let  $G_1=G_2=G-\{v,u_1,u_2\}$ and $G_3=G-\{v\}$. By the induction hypothesis we can determine $L$-choosability of $G$ in    
  $O((1.3196)^{n_3-6+.5(n_2+2)}+(1.3196)^{n_3-6+.5(n_2+3)}+(1.3196)^{n_3-5+.5(n_2+4)})$ time, which is a subset of $O(1.3196^{n_3+.5n_2})$, as desired.

If $\{c_1,c_2\}\neq \{1,2\}$, then by symmetry we may suppose $\{c_1,c_2\}=\{1,3\}$. In this case let $G_1=G-\{v,u_1,u_2\}$, $G_2=G-\{v,u_1\}$, and $G_3=G-\{v,u_2\}$. By the induction hypothesis we can determine $L$-choosability of $G$ in  
$O((1.3196)^{n_3-6+.5(n_2+3)}+(1.3196)^{n_3-6+.5(n_2+4)}+(1.3196)^{n_3-5+.5(n_2+3)})$ time, which is a subset of $O(1.3196^{n_3+.5n_2})$, as desired.

Note that $u'_1$ might be a neighbor of $v$. In this case when $v$ gets color 1 or 2, a neighbor of $v$ with list of size 3 (the vertex $u'$) reduces to a vertex with list of size 1.  As a result, when $\{c_1,c_2\}=\{1,2\}$,  we can determine $L$-choosability of $G$ in  
$O((1.3196)^{n_3-5+.5(n_2+1)}+(1.3196)^{n_3-5+.5(n_2+1)}+(1.3196)^{n_3-5+.5(n_2+4)})$ time and when $\{c_1,c_2\}\neq \{1,2\}$,   we can determine $L$-choosability of $G$ in  
$O((1.3196)^{n_3-5+.5(n_2+2)}+(1.3196)^{n_3-5+.5(n_2+1)}+(1.3196)^{n_3-5+.5(n_2+3)})$ 
time, which are both subsets of $O(1.3196^{n_3+.5n_2})$, as desired.

A similar argument as the argument in Subcase 1 applies if  $u_2$ has another neighbor with a list of size 3. Therefore for the rest of cases we may suppose the only neighbors of $u_1$ and $u_2$ with  list of size 3 is $v$.

\noindent
\textbf{Subcase 2}: Each of the vertices $u_1$ and $u_2$ has degree 2 in $G$. In this case if $v$ gets color 1 or 2, then $u_1$ must have color 2 or 1, respectively. If $v$ gets color 3, then $u_1$ keeps its list 
($\{1,2\}$). However it only needs to avoid one color that is  the color of its other neighbor. Therefore we can remove it from the graph, color the smaller graph, if possible, and then extend the coloring to a coloring of $G$ by giving $u_1$ an appropriate color in $\{1,2\}$. A similar argument applies to vertex $u_2$.

Therefore in all the three possibilities of the color of $v$, we can reduce the problem into an $L'$-list choosability of the graph $G-\{v,u_1,u_2\}$. Hence by the induction hypothesis we can determine $L$-choosability of $G$ in time 
$O(3(1.3196)^{n_3-5+.5(n_2+2)})$, which is a subset of $O(1.3196^{n_3+.5n_2})$, as desired.

Therefore for the rest of cases we may suppose $u_1$ has degree at least 3 in $G$.

\noindent
\textbf{Subcase 3}: $u_1$ has a neighbor $w_1$ with list $\{1,2\}$ and a neighbor $w_2$ with list $\{2,3\}$, and the list on $u_2$ is either $\{1,3\}$ or $\{2,3\}$ (Figure 5).

\begin{figure}[htp]
    \centering
    \includegraphics[width=7cm]{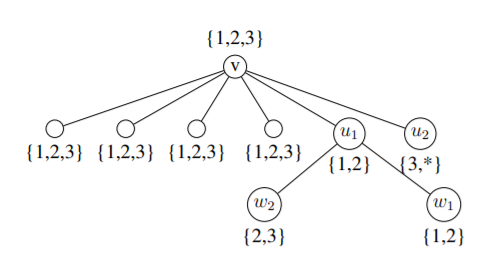}
    \caption{When $L(w_1)=\{1,2\}$ and $L(w_2)=\{2,3\}$ and $3\in L(u_2)$.}
    \label{TL5}
\end{figure}

If $v_1$ gets color 1, then we reduce the problem into a list choosability of $G-\{v, u_1,w_1,w_2\}$ (Figure 6). If $v_1$ gets color 2, then we reduce the problem into a list choosability of $G-\{v,u_1,w_1\}$, and when $v$ gets color 3, we reduce the problem into a list choosability of $G-\{v,u_2\}$.

\begin{figure}[htp]
    \centering
    \includegraphics[width=7cm]{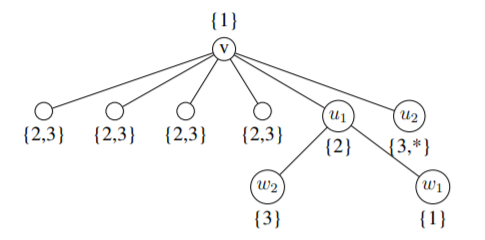}
    \caption{When $v$ gets color 1 in Subcase 3.}
    \label{TL6}
\end{figure}

 By the induction hypothesis it takes $O((1.3196)^{n_3-5+.5(n_2+1)}+(1.3196)^{n_3-5+.5(n_2+2)}+(1.3196)^{n_3-5+.5(n_2+3)})$ time to determine $L$-choosability of $G$,
which is a subset of $O(1.3196^{n_3+.5n_2})$, as desired.

\noindent
\textbf{Subcase 4}: $u_1$ has a neighbor $w_1$ with list $\{1,2\}$ and a neighbor $w_2$ with list $\{1,3\}$, and the list on $u_2$ is either $\{1,3\}$ or $\{2,3\}$. The argument in this case is very similar to that in Subcase 3.

\noindent
\textbf{Subcase 5}: $u_1$ has a neighbor $w_1$ with list $\{1,2\}$ and a neighbor $w_2$ with list $\{2,3\}$, and the list on $u_2$ is $\{1,2\}$ (Figure 7).

\begin{figure}[htp]
    \centering
    \includegraphics[width=7cm]{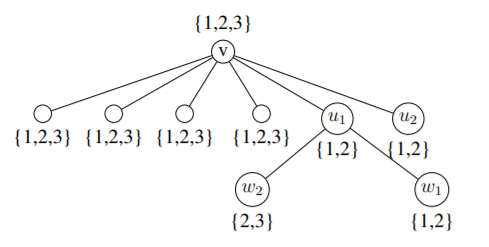}
    \caption{When $L(w_1)=L(u_2)=\{1,2\}$ and $L(w_2)=\{2,3\}$.}
    \label{TL7}
\end{figure}

If $v_1$ gets color 1 (Figure 8), then we reduce the problem into a list choosability of $G-\{v,u_1,w_1,w_2,u_2\}$. If $v_1$ gets color 2, then we reduce the problem into a list choosability of $G-\{v,u_1,w_1,u_2\}$, and when $v$ gets color 3, we reduce the problem into a list choosability of $G-\{v\}$. By the induction hypothesis it takes $O(
(1.3196)^{n_3-5+.5(n_2)}+(1.3196)^{n_3-5+.5(n_2+1)}+(1.3196)^{n_3-5+.5(n_2+4)})$ time to determine $L$-choosability of $G$. 
Since $
(1.3196)^{n_3-5+.5(n_2)}+(1.3196)^{n_3-5+.5(n_2+1)}+(1.3196)^{n_3-5+.5(n_2+4)}\leq 
1.3196^{n_3+.5n_2}$ the assertion holds in this case.

\begin{figure}[htp]
    \centering
    \includegraphics[width=7cm]{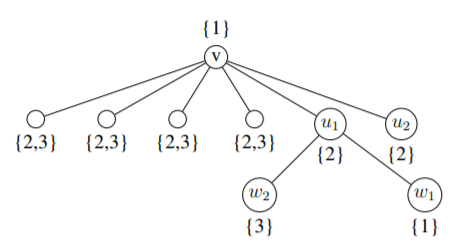}
    \caption{When $v$ gets color 1 in Subcase 5.}
    \label{TL8}
\end{figure}

\noindent
\textbf{Subcase 6}: $u_1$ has a neighbor $w_1$ with list $\{1,2\}$ and a neighbor $w_2$ with list $\{1,3\}$, and the list on $u_2$ is $\{1,2\}$. The argument in this case is very similar to that in Subcase 5.

\noindent
\textbf{Subcase 7}: $u_1$ has two neighbors $w_1$ and $w_2$ both with lists $\{1,2\}$ and the list of $u_2$ is $\{1,3\}$ or $\{2,3\}$.

If $v_1$ gets color 1 or 2, then we reduce the problem into a list choosability of $G-\{v,u_1,w_1,w_2\}$ and when $v$ gets color 3, we reduce the problem into a list choosability of $G-\{v,u_2\}$. By the induction hypothesis it takes $O(2(1.3196)^{n_3-5+.5(n_2+1)}+(1.3196)^{n_3-5+.5(n_2+3)})$ time to determine $L$-choosability of $G$,
which is a subset of $O(1.3196^{n_3+.5n_2})$, as desired.

\noindent
\textbf{Subcase 8}: $u_1$ has two neighbors $w_1$ and $w_2$ both with lists $\{1,2\}$ and the list of $u_2$ is $\{1,2\}$.

If $v_1$ gets color 1 or 2, then we reduce the problem into a list choosability of $G-\{v,u_1,w_1,w_2,u_2\}$ and when $v$ gets color 3, we reduce the problem into a list choosability of $G-\{v\}$. By the induction hypothesis it takes $O(2(1.3196)^{n_3-5+.5(n_2)}+(1.3196)^{n_3-5+.5(n_2+4)})$ time to determine $L$-choosability of $G$,
which is a subset of $O(1.3196^{n_3+.5n_2})$, as desired.

\noindent
\textbf{Subcase 9}: $u_1$ has a neighbor $w_1$ with list $\{1,3\}$ and a neighbor $w_2$ with list $\{2,3\}$, and the list on $u_2$ is either $\{1,3\}$ or $\{2,3\}$. We may suppose the list on $u_2$ is $\{1,3\}$ (Figure 9). The proof of the other case is similar. 

\begin{figure}[htp]
    \centering
    \includegraphics[width=7cm]{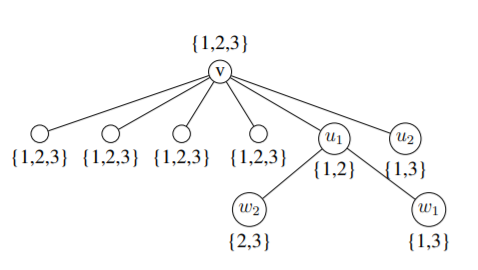}
    \caption{When $L(w_1)=\{1,2\}$ and $L(w_2)=\{1,3\}$ and $3\in L(u_2)$.}
    \label{TL9}
\end{figure}

For the cases that $v$ gets color 1, 2, or 3, we reduce the problem into choosability problems of the smaller graphs $G-\{v,u_1,w_2,u_2\}$, $G-\{v,u
_1,w_1\}$, and $G-\{v,u_2\}$, respectively. Figure 10 explains the case $v$ gets color 1.

\begin{figure}[htp]
    \centering
    \includegraphics[width=7cm]{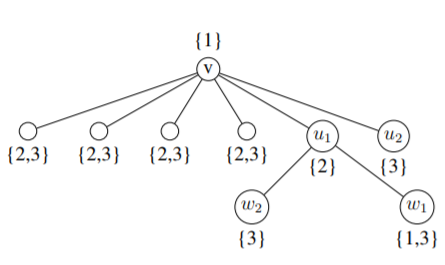}
    \caption{When $v$ gets color 1 in Subcase 9.}
    \label{TL10}
\end{figure}

 As a result the induction hypothesis implies that $L$-choosability of $G$ can be determined in time $O((1.3196)^{n_3-5+.5(n_2+1)}+(1.3196)^{n_3-5+.5(n_2+2)}+(1.3196)^{n_3-5+.5(n_2+3)})$, which is a subset of $O(1.3196^{n_3+.5n_2})$, as desired.

\noindent
\textbf{Subcase 10}: $u_1$ has a neighbor $w_1$ with list $\{1,3\}$ and a neighbor $w_2$ with list $\{2,3\}$, and the list on $u_2$ is $\{1,2\}$. 

If $u_2$ has degree 1 in $G$, then for the cases that $v$ gets color 1, 2, or 3, we reduce the problem into choosability problems of the smaller graphs $G-\{v,u_1,w_2,u_2\}$, $G-\{v,u
_1,w_1,u_2\}$, and $G-\{v,u_2\}$, respectively. As a result the induction hypothesis implies that $L$-choosability of $G$ can be determined in time $O(2(1.3196)^{n_3-5+.5(n_2+1)}+(1.3196)^{n_3-5+.5(n_2+3)})$, which is a subset of $O(1.3196^{n_3+.5n_2})$, as desired. 

If $u_2$ has degree more than 1 in $G$, then we may suppose a neighbor of $u_2$ that is different from $v$ has color $1$ on his list. As a result, for the cases that $v$ gets color 1, 2, or 3, we reduce the problem into choosability problems of the smaller graphs $G-\{v,u_1,w_2,u_2\}$, $G-\{v,u
_1,w_1,u_2\}$, and $G-\{v\}$, respectively. As a result the induction hypothesis implies that $L$-choosability of $G$ can be determined in time $O((1.3196)^{n_3-5+.5(n_2+1)}+(1.3196)^{n_3-5+.5(n_2)}+(1.3196)^{n_3-5+.5(n_2+4)})$, which is a subset of $O(1.3196^{n_3+.5n_2})$, as desired.

\noindent
\textbf{Subcase 11}: $u_1$ has two neighbors $w_1$ and $w_2$ both with lists $\{1,3\}$ or both with lists $\{2,3\}$ and the list of $u_2$ is $\{1,3\}$ or $\{2,3\}$ (Figure 11).

\begin{figure}[htp]
    \centering
    \includegraphics[width=7cm]{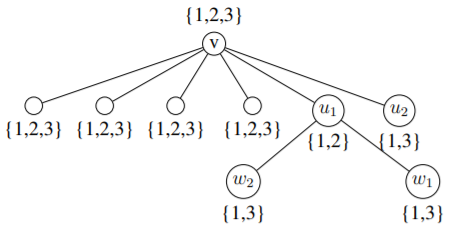}
    \caption{When $L(w_1)=L(w_2)=L(u_2)=\{1,3\}$.}
    \label{TL11}
\end{figure}

We may suppose $w_1$ and $w_2$ both have lists $\{1,3\}$, the proof of the other case is similar. We may also suppose that $u_1$ has no neighbor with list $\{2,3\}$, since otherwise the above cases can be applied.  For the cases that $v$ gets color 1, 2, or 3, we reduce the problem into choosability problems of the smaller graphs $G-\{v,u_1\}$, $G-\{v,u
_1,w_1,w_2\}$, and $G-\{v,u_1,u_2\}$, respectively. Figure 12 explains this for the case we give $v$ color $1$.

\begin{figure}[htp]
    \centering
    \includegraphics[width=7cm]{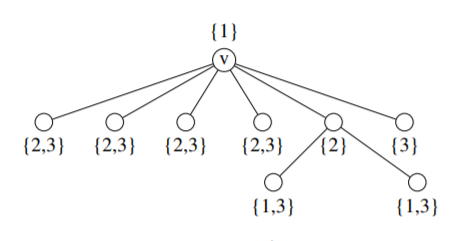}
    \caption{When $v$ gets color 1 in Subcase 11.}
    \label{TL12}
\end{figure}

 Note that for the case that $v$ gets color 3, we can give $u_1$ color 2 because none of its other neighbors' list contain color 2. By the induction hypothesis $L$-choosability of $G$ can be determined in time $O((1.3196)^{n_3-5+.5(n_2+3)}+(1.3196)^{n_3-5+.5(n_2+1)}+(1.3196)^{n_3-5+.5(n_2+2)})$, which is a subset of $O(1.3196^{n_3+.5n_2})$, as desired.

\noindent
\textbf{Subcase 12}: $u_1$ has two neighbors $w_1$ and $w_2$ both with lists $\{1,3\}$ or both with lists $\{2,3\}$ and the list of $u_2$ is $\{1,2\}$.

We may suppose $w_1$ and $w_2$ both have lists $\{1,3\}$, the proof of the other case is similar. We may also suppose that $u_1$ has no neighbor with list $\{2,3\}$, since otherwise the above cases can be applied.  For the cases that $v$ gets color 1, 2, or 3, we reduce the problem into choosability problems of the smaller graphs $G-\{v,u_1,u_2\}$, $G-\{v,u
_1,w_1,w_2,u_2\}$, and $G-\{v,u_1\}$, respectively. Note that for the case that $v$ gets color 3, we can give $u_1$ color 2 because none of its other neighbors' list contain color 2. By the induction hypothesis $L$-choosability of $G$ can be determined in time $O((1.3196)^{n_3-5+.5(n_2+2)}+(1.3196)^{n_3-5+.5(n_2)}+(1.3196)^{n_3-5+.5(n_2+3)})$, which is at most $O(1.3196^{n_3+.5n_2})$, as desired.

\noindent
\textbf{Case 4}: None of the above cases. Let $A$ be the set of vertices with lists of size 2 in $G$ and let $B$ be the vertices with lists of size 3 in $G$. Since we are supposing that none of the above cases happen, each vertex of $B$ has at least 3 neighbors in $A$.

If there is a vertex in $A$ having at least 2 neighbors in $B$, call it $a_1$. Now if there is a vertex in $A-\{a_1\}$ having at least 2 neighbors in $B-N(a_1)$, then call this vertex $a_2$. At step $i$ if there is a vertex in $A-\{a_1,\ldots,a_{i-1}\}$  having at least two neighbors in $B-\cup_{i=1}^{i-1}N(a_i)$, name the vertex $a_i$. Continue the process until no such vertex exists. Suppose the process stops at step $x_1$ and let $X_1$ be the set of these vertices. i.e. $X_1=\{a_1,\ldots,a_{x_1}\}$. 

Let $Y_1$ be the set of neighbors of $X_1$ in $B$ and let $X_2=A-X_1$, $Y_2=B-Y_1$, $x_2=|X_2|$, $y_1=|Y_1|$, and $y_2=|Y_2|$ (Figure 13). By the choice of $X_1$ and $Y_1$ we have 

$$y_1\geq 2x_1.$$

 Note that each vertex in $X_2$ has at most one neighbor in $Y_2$, since otherwise we can extend $X_1$ to a larger set. Moreover, each vertex in $Y_2$ has at least three neighbors in $X_2$, because as we argued above each vertex with a list of size 3 has at least three neighbors with lists of size 2 (Figure 13). Therefore 
 
 $$y_2\leq 3x_2.$$
 
\begin{figure}[htp]
    \centering
    \includegraphics[width=7cm]{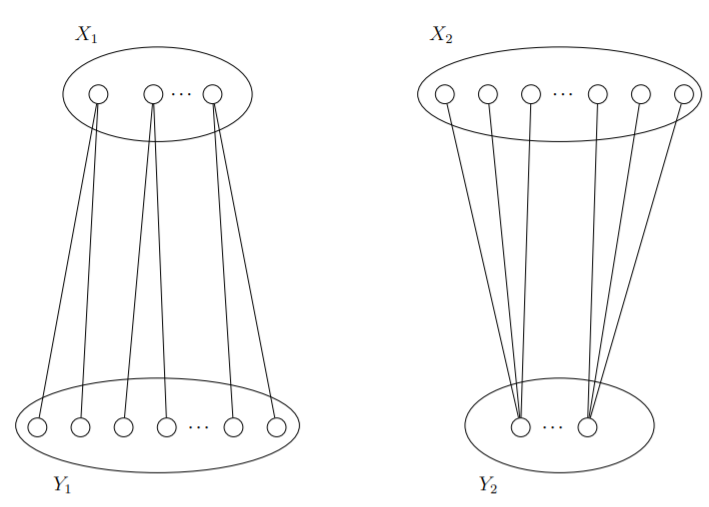}
    \caption{Partitioning $V(G)$ into vertices with lists of size 2 and size 3.}
    \label{TL13}
\end{figure}

 Now we break $L$-choosability of $G$ into list choosability of $2^{x_1+x_2}$ smaller graphs such that each vertex in those graphs has a list of size 1 or 2. In fact each vertex in $X_1$ has a list of size 2, so it has two color options. Once we fix the coloring of the vertices in $X_1$, the lists on the vertices in $Y_1$ go down to size 1 or 2. Moreover, for each vertex in $Y_2$, we consider the two possibilities that the vertex gets color 1 or another color. If the color of $v$ is not 1, then its list becomes $\{2,3\}$, which has size 2. Each of the resulting $2^{x_1+y_2}$ graphs has vertices of size at most 2. Hence by [2] we can determine in polynomial time if it has a proper coloring from its lists or not. Hence the complexity of determining if $G$ is $L$-choosable or not is $O(2^{x_1+y_2})$. Since $y_1\geq 2x_1$ and $y_2\leq 3x_2$, we have 
 
 $$2^{x_1}\leq 1.3196^{y_1+.5x_1}$$
 
 and
 
 $$2^{y_2}\leq 1.3196^{y_2+.5x_2}$$

Therefore $O(2^{x_1+y_2})\subseteq O(1.3196^{y_1+y_2+.5x_1+.5x_2})=O(1.3196^{n_3+.5n_2})$, as desired.


\begin{thebibliography}{9}
\frenchspacing


















\bibitem{BE} R. Beigel and D. Eppstein, 3-coloring in time $O(1.3289^n)$, \textit{J. Algorithms}, 54:2, 168--204, 2005.


\bibitem{BHK} A. Bj\"orklund, T. Husfeldt, and M. Koivisto, Set partitioning via inclusion–exclusion, \textit{SIAM J. Comput. 39} (2009), 546--563.


\bibitem{C} N. Christofides, An Algorithm for the Chromatic Number of a Graph,  \textit{Computer J.}, 14, 38--39, 1971.


\bibitem{CJP} N. Crawford, S.Jahanbekam, K.Potika, Improved algorithm to determine 3-colorability of graphs with minimum degree at least 7,  Discrete Applied Mathematics, to appear.

\bibitem{MTB}  R. E. Miller, J. W. Thatcher, and J.D. Bohlinger,  Complexity of Computer Computations,  \textit{New York: Plenum.} pp. 85–103.



\bibitem{KT} J. Kratochv\'ila and Z. Tuza, Algorithmic complexity of list colorings,\textit{ Discrete Applied Mathematics}, 
Volume 50, Issue 3,  297--302, 1994.












\bibitem{W} 
D. B. West,\textit{ Introduction to Graph Theory}, Second edition,  Published by Prentice Hall 1996, 2001. ISBN 0-13-014400-2.

\end{thebibliography}
\end{document}